\documentclass[11pt]{article}

\usepackage{amssymb,amsmath,amsfonts,amssymb,numbysec}
\usepackage{graphics,graphicx,color,epsfig}

\hsize \textwidth
\advance \hsize by -\marginparwidth
\evensidemargin \oddsidemargin
\usepackage{amssymb}
\advance\hoffset by 5mm

\def\@abssec#1{\vspace{.05in}\footnotesize \parindent .2in
{\bf #1. }\ignorespaces}

\newtheorem{theorem}{Theorem}[section]
\newtheorem{lemma}[theorem]{Lemma}
\newtheorem{proposition}[theorem]{Proposition}

\def \Rm {\mathbb R}

\newcommand{\eps}{\varepsilon}
\newcommand{\E}{\mathbb E}

\newcommand{\pdr}[2]{\dfrac{\partial{#1}}{\partial{#2}}}

\newcommand{\vz}{\mathbf z}

\newcommand{\bk}{k}\newcommand{\vy}{y} \newcommand{\vx}{x}
\newcommand{\bx}{x}  \newcommand{\bp}{p}

\newcommand{\commentout}[1]{}

\renewcommand{\thefootnote}{{\arabic{footnote}}}

 \renewcommand{\arraystretch}{1.5}
\numberbysection
\title{Wave field correlations in weakly mismatched random media }

\author{Guillaume Bal \footnote{Department of
   Applied Physics and Applied Mathematics, Columbia University,
   New York NY, 10027; gb2030@columbia.edu}
\and Lenya Ryzhik \footnote{Department of Mathematics, University
of Chicago, Chicago IL, 60637; ryzhik@math.uchicago.edu} }

\begin{document}

\maketitle


\begin{abstract}
  This paper concerns the derivation of a Fokker-Planck equation for
  the correlation of two high frequency wave fields propagating in two
  different random media. The mismatch between the random media need
  be small, on the order of the wavelength, and their correlation
  length need be large relative to the wavelength. The loss of
  correlation caused by the mismatch in the random media is quantified
  and the limit process for the phase difference is obtained. The
  derivation is based on a random Liouville equation to model high
  frequency correlations and on the method of characteristics to
  characterize mixing in the random Liouville equation. Applications
  of such correlation loss include the monitoring in time of random
  media and the analysis of time reversed waves in changing
  heterogeneous domains.
\end{abstract}

\renewcommand{\thefootnote}{\fnsymbol{footnote}}
\renewcommand{\thefootnote}{\arabic{footnote}}

\renewcommand{\arraystretch}{1.1}





\section{Introduction}
\label{sec:intro}

The energy density of high frequency waves propagating in highly
heterogeneous media can be modeled by a Fokker-Planck equation in the
phase space, i.e., the space of positions and momenta. We refer to
e.g. \cite{BKR-liouv} for a mathematical derivation of such a
macroscopic model for the wave energy density. The main assumption on
the heterogeneous medium is that it is a random medium with a very
large correlation length relative to the typical wavelength in the
system.  The Fokker-Planck equation may be seen as a
highly-peaked-forward-scattering approximation to the radiative
transfer equations, which are also used in the modeling of the energy
density of waves in heterogeneous media when correlation length and
wavelength are comparable; see e.g.
\cite{B-WM-05,ishimaru-93,RPK-WM}.

Such macroscopic models for waves in random media can more generally
be used to quantify the correlation function of two wave fields
propagating in possibly two different media \cite{B-WM-05}. This has
applications in the temporal monitoring of the statistical properties
of random media as well as in the analysis of the refocusing of time
reversed waves \cite{BV-MMS-04,BR-DCDS-A05,LVKBC-TRCM-06}. This paper
analyzes the effect of changes in the random media on the two-field
correlation. In the Fokker-Planck regime, the two-field correlation
decays as the two media separate and we present a quantitative
estimate of such a decay.

Although the results in this paper generalize to fairly large classes
of waves such as e.g., acoustic waves as in \cite{BKR-liouv}.
electromagnetic waves, and elastic waves (see \cite{RPK-WM}) we
restrict ourselves to the case of a scalar Schr\"odinger equation to
simplify.  This models the effect of heterogeneities on a single
particle represented by a quantum wave function. The changes in the
heterogeneous medium considered here are sufficiently small so that
the correlation of the two fields, one propagating in the unperturbed
medium and the other one propagation in the perturbed medium, still
satisfies a random Liouville equation in the high frequency
limit. However, solutions of the Liouville equation are no longer real
and have a complex phase. Their phase is responsible for the decoherence
of the two wave fields, and is driven by the media mismatch. The
evolution of the correlation function is described in terms of the
particle position $X^\delta(t)$, momentum $K^\delta(t)$ and phase
difference $Z^\delta(t)$. The techniques used in \cite{BKR-liouv} are
then applied to this new random Liouville equation to study the limit
of the joint process $(X^\delta(t),K^\delta(t),Z^\delta(t))$. The
limiting equation, of Fokker-Planck type, is obtained in the vanishing
limit of the correlation length $\delta$ of the heterogeneous
medium. Its derivation is based on the mixing properties of the
bi-characteristics of a random Hamiltonian, as in
\cite{BKR-liouv}, and on the mixing properties of a highly oscillatory
functional of such bi-characteristics, which is the main new result
obtained in this paper. In particular, we show that the phase difference
$Z^\delta(t)$ converges to a Brownian motion on the real line.

The rest of the paper is organized as follows. Section \ref{sec:rmLe}
presents the random Schr\"odinger equations and the limiting random
Liouville equation for the Wigner transform of the two wave fields,
which is the Fourier transform in the difference variable of the
correlation function of the two fields. We also review in Section
\ref{sub:review} the results obtained in \cite{BKR-liouv} while adapting
them to the case of a random Schr\"odinger equation. Finally, in
Section \ref{sub:formal} we derive the generalized Fokker-Planck
equations (see \eqref{formal-W} below) by formal analysis. Section
\ref{sec:liouv} presents Theorem \ref{thm-new}, the main result of
this paper that describes the limit of the joint process
$(X^\delta(t),K^\delta(t),Z^\delta(t))$. The Fokker-Planck equation is
deduced from the limiting law (as the correlation length of the random
medium goes to zero).  The proof of Theorem \ref{thm-new} is presented
in Section \ref{sec:proof}.

{\bf Acknowledgment.} LR and GB have been partially supported by the
ONR and Alfred P. Sloan Fellowships. GB has been also supported by NSF
grant DMS-0239097.

\section{The random momenta Liouville equation}
\label{sec:rmLe}

We consider the evolution of the correlation function of two
solutions of the Schr\"odinger equation with a small mismatch
between the random potentials but with the same initial data. The
function $\psi_\eps$ satisfies
\begin{equation}\label{schr-psi}
i\eps\pdr{\psi_\eps}{t}+\frac{\eps^2}{2}\Delta\psi_\eps-V_\delta(x)\psi_\eps=0
\end{equation}
and the function $\phi_\eps$ satisfies
\begin{equation}\label{schr-phi}
i\eps\pdr{\phi_\eps}{t}+\frac{\eps^2}{2}\Delta\phi_\eps-[V_\delta(x)+\eps
S_\delta(x)]\phi_\eps=0.
\end{equation}
Both of the functions $\phi_\eps$ and $\psi_\eps$ satisfy
initially
\begin{equation}\label{phi-indata}
\psi_\eps(0,x)=\phi_\eps(0,x)=\phi_0^\eps(x).
\end{equation}
The family $\phi_\eps^0$ is $\eps$-oscillatory and compact at
infinity \cite{GMMP}. The random potentials $V_\delta$ and $S_\delta$ vary on
a scale $\delta$ that is much larger than the wave length $\eps$
of the initial data but is much smaller than the overall
propagation distance that is of the order $O(1)$:
$\eps\ll\delta\ll 1$. To keep a non-trivial correlation of
$\psi_\eps$ and $\phi_\eps$ the mismatch of the potentials has to
be weak  -- hence the coefficient $\eps$ in front of $S_\delta$.
We will see that in order to produce an order one contribution we
will have eventually to take
$S_\delta(x)=\delta^{-1/2}S(x/\delta)$ making the overall strength
of the mismatch be of the order $O(\eps/\sqrt{\delta})$.

In order to study the correlation of $\psi_\eps$ and $\phi_\eps$
we introduce the cross Wigner transform as
\begin{equation}\label{def-wigner}
W_\eps(t,\vx,\bk)=\int e^{i\bk\cdot\vy }
\psi_\eps\left(t,\vx-\frac{\eps\vy}{2}\right)
\bar\phi_\eps\left(t,\vx+\frac{\eps\vy}{2}\right)
\frac{d\vy}{(2\pi)^{d}}.
\end{equation}
When $\phi_\eps=\psi_\eps$ then $W_\eps$ is real and is usually
interpreted as the phase space energy density of the family
$\phi_\eps(x)$. However, the distribution $W_\eps(t,x,k)$ does not
have to be real if $\phi_\eps\neq \psi_\eps$. Its phase measures the
decoherence of the functions $\phi_\eps$ and $\psi_\eps$. The basic
properties of the Wigner transforms can be found in \cite{GMMP,LP}.

In order to obtain an equation for $W_\eps$ we differentiate the
Wigner transform with respect to time:
\begin{eqnarray*}
&&\pdr{W_\eps}{t}+\bk\cdot\nabla_\vx W_\eps\\&&=\frac{1}{i\eps} \int
e^{i\bk\cdot\vy} \left[V_\delta\left(x-\frac{\eps
y}{2}\right)-V_\delta\left(x+\frac{\eps y}{2}\right)-\eps
S_\delta\left(x+\frac{\eps y}{2}\right)\right]\psi_\eps\left(x-\frac{\eps
y}{2}\right)
\bar\phi_\eps\left(x+\frac{\eps y}{2}\right) \frac{d\vy}{(2\pi)^{d}}.
\end{eqnarray*}
Passing to the limit $\eps\to 0$ we obtain an equation for the
distribution $W_\delta(t,x,k)$, the weak limit of $W_\eps$ as
$\eps\to 0$:
\begin{equation}\label{delta-liouv}
\pdr{W_\delta}{t}+\bk\cdot\nabla_\vx W_\delta-\nabla
V_\delta(\vx)\cdot\nabla_\bk W_\delta= iS_\delta(\vx)W_\delta.
\end{equation}
In order to obtain a non-trivial limit of $W_\delta(t,x,k)$ as the
correlation length $\delta\to 0$ we choose the random potential
$V_\delta(x)$ and the mismatch $S_\delta(x)$ to be of the form
\[
V_\delta(\vx)= {\sqrt{\delta}}V\left(\frac{\vx}{\delta}\right),~~~
S_\delta(\vx)=\frac{1}{\sqrt{\delta}}S\left(\frac{\vx}{\delta}\right).
\]
Then (\ref{delta-liouv}) becomes
\begin{eqnarray}\label{eq-liouv}
&&\pdr{W_\delta}{t}+\bk\cdot\nabla_\vx
W_\delta-\frac{1}{\sqrt{\delta}}\nabla
V\left(\frac{\vx}{\delta}\right)\cdot\nabla_\bk W_\delta=
\frac{i}{\sqrt{\delta}}S\left(\frac{\vx}{\delta}\right)W_\delta\\
&&W(0,x,k)=W_0(x,k).\nonumber
\end{eqnarray}
The initial data $W_0(x,k)$ is simply the limit Wigner measure of the
family $\phi_\eps^0$. Equation (\ref{eq-liouv}) is the starting point
of our analysis. We note that if we consider the initial data for the
Schr\"odinger equation as a mixture of states then the error bound for
the approximation for the Wigner transform $W_\eps$ by the solution of
(\ref{eq-liouv}) may be estimated and the sequential limits $\eps\to
0$ first, and $\delta\to 0$ second may be replaced by a joint limit
$(\eps,\delta)\to 0$; see \cite{BKR-liouv,LP} for details. This may be
done in a certain region of the $(\eps,\delta)$-plane that ensures
that the scale separation $\eps\ll\delta\ll 1$ is kept under control.
We will not pursue this avenue in this paper in order to avoid
unnecessary technical complications.

\subsection{A review of the case in the absence of a mismatch}
\label{sub:review}

We first recall the known results of \cite{Kesten-Papanico} (see also
\cite{KR-Diff} for a longer time scale analysis) when there is no
potential mismatch, that is, when $S=0$ and $d\ge 3$. We restrict our
attention in this paper also to the case $d\ge 3$ though a
generalization using the results of \cite{DGL-CMP87} and \cite{KR-2d}
is possible. If $S=0$ then solution of (\ref{eq-liouv}) is given
explicitly in terms of the random characteristics
$(X^\delta(t),K^\delta(t))$: define
\begin{equation}\label{liouv-char}
\frac{dX^\delta(t)}{dt}=-K^\delta(t),~~~
\frac{dK^\delta(t)}{dt}=\frac{1}{\sqrt{\delta}}\nabla
V\left(\frac{X^\delta(t)}{\delta}\right),~~X^\delta(0)=x,~K^\delta(0)=k.
\end{equation}
We have $W_\delta(t,x,k)=W_0(X^\delta(t),K^\delta(t))$. The
characteristic trajectories satisfy a limit theorem; the process
$K^\delta(t)$ converges to a Brownian motion $K(t)$ on the sphere
$\{|k|=|k(0)|\}$ and $X^\delta(t)$ converges in law to
\[
X(t)=x+\int_0^tK(s)ds.
\]
Before formulating the limit theorem we describe first the necessary
assumptions on the random potentials $V$ and $S$; we will not, of
course, need the assumptions on $S$ in the present section but we will
need them later on and it is convenient to put them here.

Let $(\Omega,\Sigma,\mathbb P)$ be a probability space, and let $\E$
denote the expectation with respect to $\mathbb P$.  Let us denote the
joint process $F(x,\omega)=(V(x,\omega),S(x,\omega))$ and assume that
the random field $F$ is stationary in the first variable. This means
that for any shift $x\in\Rm^d$ and a collection of points
$x_1,\ldots,x_n\in\Rm^d$ the laws of $(F(x_1+ x), \ldots,F( x_n+ x )
)$ and $(F(x_1), \ldots,F(x_n) )$ are identical. In addition, we
assume that $\mathbb E \{S(x)\}=\E \{V(x)\}=0$ for all $x\in\Rm^d$,
the realizations of $V(x)$ and $S(x)$ are $\mathbb P$ a.s.
$C^2$-smooth in $x$ and they satisfy
\begin{equation}\label{d-i}
 M:=\max\limits_{|\alpha|\le 2}\,\mathop{\mbox{ess-sup}}
\limits_{( x, \omega)\in \Rm^d\times\Omega} |\partial_x^\alpha
F(x,\omega)|<+\infty.
\end{equation}

We suppose further that the random field $F(x,\omega)$ is strongly
mixing in the uniform sense.  More precisely, for any $R>0$ we let
${\cal C}_{R}^i$ and ${\cal C}_{R}^e$ be the $\sigma$-algebras
generated by random variables $F(x)$ for all $ x\in \mathbb B_R$ and $
x\in \mathbb B_R^c$ respectively. The uniform mixing coefficient
between the $\sigma$-algebras is
\begin{equation}\label{mix-coeff}
\phi(\rho):=\sup[\,|\mathbb P(B)-\mathbb P(B|A)|:\,R>0,\,A\in
{\cal C}_{R}^i,\,B\in {\cal C}_{R+\rho}^e\,],
\end{equation}
for all $\rho>0$. We suppose that $\phi(\rho)$ decays faster than
any power: for each $p>0$
\begin{equation}\label{DR}
h_p:=\sup\limits_{\rho\ge0}\rho^p\phi(\rho)<+\infty.
\end{equation}
The two-point spatial correlation $2\times 2$ tensor of the random
field $F$ is denoted by $R(y)$ and has components
\[
R^{VV}(y)=\E[V(0)V(y)],~ R^{VS}(y)=\E\{V(0)S(y)\},~
R^{SV}(y)=\E\{S(0)V(y)\},~R^{SS}(y)=\E\{S(0)S(y)\}.
\]
Note  that (\ref{DR}) implies that for each $p>0$
\begin{equation}\label{53102-intro}
h_{p} :=\, \sum\limits_{|\alpha|\le 4}\sup\limits_{y \in\Rm^d
}(1+| y|^2)^{p/2} |\partial_y^\alpha
 R( y)|<+\infty.
\end{equation}
We also assume that the correlation tensor $R(y)$ is of the
$C^\infty$-class  and that
\begin{equation}\label{non-trivial}
\hbox{$\hat R^{VV}(k)$ does not vanish identically on any
hyperplane $H_\bp=\{\bk:~(\bk\cdot\bp)=0\}$.}
\end{equation}
Here $\hat R^{VV}(\bk)=\int R^{VV}(\vx)\exp(-i k\cdot x)d x$ is
the power spectrum of $V$.

Let us define the diffusion matrix $D_{mn}$ by
\begin{equation}\label{diff-matrix1-main}
D_{mn}(k)=-\frac{1}{2}
\int_{-\infty}^\infty\frac{\partial^2R^{VV}(ks)}{\partial
x_n\partial x_m} ds=-\frac{1}{2|k|}\int_{-\infty}^\infty
\frac{\partial^2R^{VV}(s\hat{k} )}{\partial x_n\partial x_m}ds
,\quad\,m,n=1,\ldots,d,~~~\hat{k} =k/|k|.
\end{equation}
Then we have the following result.
\begin{theorem}\label{main-thm}\cite{Kesten-Papanico}
Let $W^\delta$ be the solution of $(\ref{eq-liouv})$ with the
initial data $W_0(x,k)$ supported in a compact set away from
$k=0$: $\hbox{supp} \{W_0\}\subset {\cal S}\times A(M)$ with
$A(M)=\{M^{-1}\le |k|\le M\}$ for some $M>0$ and a compact set
${\cal S}\subset\Rm^d$. Let the function $\bar\phi$ satisfy
\begin{eqnarray}\label{eq-mainthm-1}
&&\pdr{\bar\phi}{t}+k\cdot\nabla_x\bar\phi=
 \pdr{}{k_m}\left(D_{mn}(
k)\pdr{\bar\phi}{k_n}\right)\\
&&\bar\phi(0,x,k)=W_0(x,k).\nonumber
\end{eqnarray}
Then, there exist two constants $C(T)$ and $\alpha_0>0$ such that
\begin{equation}\label{decorrel-main1}
\sup\limits_{(t,x,k)\in [0,T]\times K} \left|\E
W^{\delta}\left(t,x,\bk\right)-
 \bar \phi(t,\vx,\bk)\right|\le
C(T)(1+\|W_0\|_{C^4})\delta^{\alpha_0}
\end{equation}
for all compact sets $K\subset{\cal A}(M)=\Rm^d\times A(M)$.
\end{theorem}
Note that
\[
  D_{nm}( k)\hat k_m= - \frac{1}{2|k|}\int_{-\infty}^\infty
\frac{\partial^2R^{VV}(s \hat k)}{\partial x_n\partial x_m}\hat
k_mds
 =
- \frac{1}{2|k|}\int_{-\infty}^\infty\frac{d}{ds}
\left(\pdr{R^{VV}(s k)}{x_n}\right)ds=0
\]
and thus the $K$-process generated by (\ref{eq-mainthm-1}) is
indeed a diffusion process on a sphere $k=\hbox{const}$, or,
equivalently, equations (\ref{eq-mainthm-1}) for different values
of $|k|$ are decoupled. It is easy to check that assumption
(\ref{non-trivial}) implies that the matrix $D(k)$ has rank $d-1$
for each $k\in{\Rm}^{d}\{0\}$. It can  be also shown that then
 equation (\ref{eq-mainthm-1}) is hypoelliptic
on the manifold $\Rm^{d}\times\{|k|=k_0\}$ for each $k_0>0$.

\subsection{A formal analysis of the momenta Liouville equation}
\label{sub:formal}

We first present a non-rigorous formal multiple scales analysis of
(\ref{eq-liouv}), which provides a short and relatively quick way to
the correct limit. We introduce a multiple scales expansion
\[
W_\delta=W(t,x,k)+\sqrt{\delta}W_1(t,x,y,k)+\delta
W_2(t,x,y,k)+\dots, ~~~y= x/\delta
\]
and insert it into (\ref{eq-liouv}). As usual we make an
additional assumption that the leading order term $W(t,x,k)$ is
deterministic and does not depend on the fast scale variable $y$.
In the leading order we obtain
\[
k\cdot\nabla_y W_1+\theta W_1=\nabla V\left({y}\right)\cdot\nabla_
k W+ {i}S\left(y\right)W.
\]
Here $\theta>0$ is an auxiliary regularizing parameter that we
will send to zero at the end. Define the correctors $\chi_j$ and
$\eta$ as mean-zero solutions of
\begin{eqnarray*}
&& k\cdot\nabla_y\chi_j+\theta\chi_j=\pdr{V}{y_j}\\
&& k\cdot\nabla_y\eta+\theta\eta=S(y).
\end{eqnarray*}
They are given explicitly by
\begin{equation}\label{b-chij}
\chi_j(y,k)=\int_0^\infty\pdr{V(y-sk)}{y_j}e^{-\theta s}ds
\end{equation}
and
\begin{equation}\label{b-eta}
\eta(y,k)=\int_0^\infty e^{-\theta s}S(y-s k)ds.
\end{equation}
The function $W_1$ is given in terms of the correctors as
\[
W_1(t,x,y,k)=\sum_{j=1}^d
\chi_j(y,k)\pdr{W(t,x,k)}{k_j}+i\eta(y,k)W(t,x,k).
\]
The equation for $W_2$ is
\[
\pdr{W}{t}+k\cdot\nabla_x W+k\cdot\nabla_y W_2=\nabla
V(y)\cdot\nabla_k W_1+iS(y)W_1.
\]
Averaging under the assumption that $\E\{ k\cdot\nabla_yW_2\}=0$
we obtain the following closed equation for the leading order term
$W$:
\begin{equation}\label{order-2}
\pdr{W}{t}+k\cdot\nabla_x W=\E\left\{\nabla V(y)\cdot\nabla_k
W_1+iS(y)W_1\right\}=J_I+J_{II}.
\end{equation}
The two terms on the right side are computed using the explicit
expressions (\ref{b-chij}) and (\ref{b-eta}) for the correctors.
The first term may be split as
\[
J_I=\E\left\{\nabla V(y)\cdot\nabla_k W_1\right\}=J_{I}^1+J_{I}^2
\]
with
\begin{eqnarray*}
&&J_{I}^1=\E\left\{\pdr{V}{y_j}(y)\pdr{}{k_j}\left[\chi_m(y, k)
\pdr{W(t, x, k)}{k_m}\right]\right\}\\
&&=\pdr{}{k_j}\left[\E\left\{\pdr{V}{y_j}(y)\int_0^\infty
\pdr{V(y-s k)}{y_m}e^{-\theta s}ds\right\}\pdr{W(t, x,
k)}{k_m}\right]=\pdr{}{k_j} \left(D_{jm}(
k)\pdr{W(t,x,k)}{k_m}\right)
\end{eqnarray*}
where the diffusion matrix $D_{jm}$ is given by
(\ref{diff-matrix1-main}). The term $J_{I}^2$ is
\begin{eqnarray*}
&&J_{I}^2=\E\left\{\pdr{V}{y_j}(y)\pdr{}{k_j}\left[i\eta(y,k)
{W(t,x,k)}\right]
\right\}=i\pdr{}{k_j}\left[\E\left\{\pdr{V}{y_j}(y)\int_0^\infty
S(y-sk) e^{-\theta s}ds\right\}{W(t,x,k)}\right]\\
&&=i\pdr{}{k_j}\left(E_j'(k)W(t,x,k)\right)
\end{eqnarray*}
with the drift
\[
E_j'(k)=\int_0^\infty\pdr{R^{SV}(sk)}{x_j}ds.
\]
Now we look at the second term in the right side of
(\ref{order-2})
\begin{equation}\label{JII}
J_{II}=\E\left\{iS(y) W_1\right\}=J_{II}^1+J_{II}^2
\end{equation}
with
\begin{eqnarray*}
&&J_{II}^1=\E\left\{iS(y)\chi_m(y, k)\pdr{W(t,x,k)}{k_m}\right\}
=i\E\left\{S(y)\int_0^\infty\pdr{V(y-s k)}{y_m}e^{-\theta
s}ds\right\}\pdr{W(t, x, k)}{k_m}\\
&&~~~~=i{E_m''}( k)\pdr{W(t, x, k)}{k_m}
\end{eqnarray*}
with
\[
E_m''=-\int_0^\infty\pdr{R^{VS}(s k)}{x_j}ds=\int_0^\infty
\pdr{R^{SV}(-s k)}{x_j}ds=\int_{-\infty}^0\pdr{R^{SV}(s
k)}{x_j}ds.
\]
Note that
\begin{eqnarray*}
&&J_I^2+J_{II}^1=i\pdr{}{k_j}\left(E_j'(\bk)W(t,\vx,\bk)\right)+iF_m(\bk)
\pdr{W(t,\vx,\bk)}{k_m}\\
&&~~~~~~~~~~~=
i(E_j'+E_j'')\pdr{W(t,\vx,\bk)}{k_j}+i(\nabla_\bk\cdot
E')W(t,\vx,\bk)=E_j\pdr{W(t,\vx,\bk)}{k_j}+FW(t,\vx,\bk)
\end{eqnarray*}
with
\begin{equation}\label{Ej-def}
E_j=E_j'+E_j''=\int_{-\infty}^\infty\pdr{R^{SV}(s k)}{x_j}ds
\end{equation}
and
\begin{equation}\label{F}
F=\nabla_k\cdot E'=\int_0^\infty s\Delta R^{SV}(s\bk)ds.
\end{equation}
The last term in (\ref{JII}) is
\[
J_{II}^2=\E\left\{iS(\vz)\eta(\vz,\bk) {W(t,\vx,\bk)} \right\}=
-\kappa(\bk)W(t,\vx,\bk)
\]
with the absorption coefficient
\begin{equation}\label{kappa-def}
\kappa(\bk)=\int_0^\infty R^{SS}(s\bk)ds.
\end{equation}
Putting together all the terms above we get the equation for $W$:
\begin{equation}\label{formal-W}
\pdr{W}{t}+k\cdot\nabla_x W=iE(k)\cdot\nabla_kW+iF(k)W+
\pdr{}{k_m}\left(D_{mn}(k)\pdr{W}{k_n}\right)-\kappa(k)W.
\end{equation}
If $S$ and $V$ are independent then $F=E=0$ and this simplifies to
\begin{equation}\label{formal-W0}
\pdr{W}{t}+k\cdot\nabla_x W=
\pdr{}{k_m}\left(D_{mn}(k)\pdr{W}{k_n}\right)-\kappa(k)W.
\end{equation}
Unfortunately, the asymptotic expansion described in this section may not
be justified rigorously.  In the next section we present a rigorous
derivation of the limit equation (\ref{formal-W}).

\section{The Liouville equation and the phase diffusion}
\label{sec:liouv}

\subsection{An example: the decorrelated case}

The purpose of this section is to obtain the limit equation
(\ref{formal-W}) as the limit of (\ref{eq-liouv}). Recall that
solution of (\ref{eq-liouv})
\begin{equation}\label{dyn0}
\pdr{W_\delta}{t}+\bk\cdot\nabla_\vx
W_\delta-\frac{1}{\sqrt{\delta}}\nabla
V\left(\frac{\vx}{\delta}\right)\cdot\nabla_\bk W_\delta=
\frac{i}{\sqrt{\delta}}S\left(\frac{\vx}{\delta}\right)W_\delta
\end{equation}
may be obtained by the method of characteristics. Along the
trajectories we have
\begin{equation}\label{dyn1}
\frac{dX^\delta}{dt}=-K^\delta,~~\frac{dK^\delta}{dt}=
\frac{1}{\sqrt{\delta}}\nabla
V\left(\frac{X^\delta}{\delta}\right),~~\frac{dZ^\delta}{dt}=
\frac{1}{\sqrt{\delta}} S\left(\frac{X^\delta}{\delta}\right)
\end{equation}
with the initial data
\[
X^\delta(0)=x,~~K^\delta(0)=k,~~Z^\delta(0)=0.
\]
Then the solution of (\ref{dyn0}) is given by
\[
W^\delta(t,x,k)=e^{iZ^\delta(t)}W_0(X^\delta(t),K^\delta(t)).
\]
The results described in the previous section tell us that
$K^\delta$ converges to a Brownian motion $K(t)$ on the sphere,
with the diffusion matrix $D_{mn}(k)$, and $X^\delta$ converges to
its time integral:
\[
K^\delta(t)\to K(t),~~X^\delta(t)\to X(t)=\vx+\int_0^tK(s)ds.
\]
We will show below that $Z^\delta$ converges to a Brownian motion
with the diffusion coefficient $\kappa( k)$ given by
(\ref{kappa-def}). In the simplest case when $V$ and $S$ are
uncorrelated  the Kolmogorov equation for the limit process
$(X(t),K(t),Z(t))$ is
\[
\pdr{f}{t}+ k\cdot\nabla_ x
f=\pdr{}{k_m}\left(D_{mn}(k)\pdr{f}{k_n}\right)
+\kappa(k)\frac{\partial^2 f}{\partial z^2}.
\]
Recall that
\[
\overline W^\delta(t,x,k):=\E\{W^\delta(t,x,k)\}=\E_{x,k,z=0}
\left(e^{iZ^{\delta}(t)}W_0(X^\delta(t),K^\delta(t))\right).
\]
Therefore, in the limit  $\delta\to 0$, the function $ \overline
W^\delta(t,x,k)$ converges to
\[
W(t,x,k)=g(t,x,k,z=0).
\]
Here the function $g$ satisfies the Kolmogorov equation
\begin{eqnarray}\label{g-eq}
&&\pdr{g}{t}+ k\cdot\nabla_x
g=\pdr{}{k_m}\left(D_{mn}(k)\pdr{g}{k_n}\right)
+\kappa(k)\frac{\partial^2 g}{\partial z^2}\\
&&g(0,x,k,z)=e^{iz}W_0(x,k).\label{g-indata}
\end{eqnarray}
It may be written as $g(t,x,k,z)=e^{iz}q(t,x,k)$, where the
function $q$ satisfies
\begin{eqnarray}\label{q-eq}
&&\pdr{q}{t}+k\cdot\nabla_x
q=\pdr{}{k_m}\left(D_{mn}(k)\pdr{q}{k_n}\right)
-\kappa(k)q\\
&&q(0,\vx,k)=W_0(\vx,\bk).\nonumber
\end{eqnarray}
We see that actually $W(t,x,k)=g(t,x,k,z=0)=q(t,x,k)$. Note that
(\ref{q-eq}) is nothing but (\ref{formal-W0}); this relates the
approach of the present section to the formal result of the previous
section. However, the Kolmogorov equation (\ref{g-eq}) provides the
description of the whole limit process $Z(t)$ while (\ref{formal-W0})
is just one of its reductions.

\subsection{The main result}

We will prove the following theorem.
\begin{theorem}\label{thm-new} The joint process
$(X^\delta(t),K^\delta(t),Z^\delta(t))$ converges in law in the
limit $\delta\to 0$ to the diffusion process $(X(t),K(t),Z(t))$
with the joint generator
\begin{equation}\label{L-gen1}
{\cal L}\phi=D_{mn}(k)\frac{\partial^2\phi}{\partial k_m\partial
k_n}+\left[D_m(k)+D_m(-k)\right]\frac{\partial^2\phi}{\partial
k_m\partial z}+D(k)\frac{\partial^2\phi}{\partial
z^2}+E_m(k)\pdr{\phi}{k_m}+E(k)\pdr{\phi}{z}-k\cdot\nabla_x\phi
\end{equation}
with the coefficients
\begin{equation}\label{Dmn1}
D_{mn}(k)=-\int\limits_{0}^\infty
 R_{mn}^{VV}\left(sk) \right)ds,~~~~D(k)=\int_0^\infty R^{SS}(sk)ds,
\end{equation}
\begin{equation}\label{Dm1}
D_m(k)=\int\limits_{0}^\infty R_m^{SV}\left(sk \right) ds,~~~~
E_m(k)=-\int_0^\infty s   \Delta R_{m}^{VV}\left(sk \right) ds,
\end{equation}
and
\begin{equation}\label{E1}
E(k)=\int_0^\infty s \Delta R^{SV}\left(sk \right)ds.
\end{equation}
\end{theorem}
The generator may be written slightly more compactly as
\begin{equation}\label{L-comp}
{\cal L}\phi=\pdr{}{k_n}\left(D_{mn}(k)\pdr{\phi}{k_m}\right)
+\pdr{}{k_m}\left(D_m(k)\pdr{\phi}{z}\right)
+\pdr{}{z}\left(D_m(-k)\pdr{\phi}{k_m}\right)+D(k)\frac{\partial^2\phi}{\partial
z^2} -k\cdot\nabla_x\phi.
\end{equation}

\subsection{A formal computation of the limit}

We first present a formal computation, which is similar to the
derivation in section \ref{sub:formal} and leads to the generator
(\ref{L-gen1}). We start with the Liouville equation including the
phase variable
\begin{equation}\label{z-eq}
\pdr{\phi}{t}+k\cdot\nabla_x\phi- \frac{1}{\sqrt{\delta}}\nabla
V\left(\frac{x}{\delta}\right)\cdot\nabla_k\phi
-\frac{1}{\sqrt{\delta}}S\left(\frac{x}{\delta}\right)\pdr{\phi}{z}=0
\end{equation}
and consider an asymptotic expansion of the form
\[
\phi(t,x,k,z)=\bar\phi(t,x,z,k)+\sqrt{\delta}\phi_1(t,x,y,z,k)
+\delta\phi_2(t,x,y,z,k)+\dots,
~~y=x/\delta.
\]
In the leading order we get:
\[
k\cdot\nabla_y\phi_1+\theta\phi_1=\nabla
V\left(y\right)\cdot\nabla_k\bar\phi +
S\left(y\right)\pdr{\bar\phi}{z}.
\]
As before $\theta>0$ is a regularizing parameter that we will send
to zero later. The leading order term $\bar\phi(t,x,k)$ is assumed
to be deterministic and independent of the fast variable $y$.
Using the correctors $\chi_j$ and $\eta$ given by (\ref{b-chij})
and (\ref{b-eta}), respectively,
we obtain an expression for $\phi_1$ as
\[
\phi_1(t,x,y,k,z)=
\chi_j(y,k)\pdr{\bar\phi}{k_j}+\eta(y,k)\pdr{\bar\phi}{z}.
\]
The equation for $\bar\phi$ is obtained in the same way as before, as
a formal solvability condition for $\phi_2$; it reads
\begin{equation}\label{interm-barphi}
\pdr{\bar\phi}{t}+k\cdot\nabla_x\bar\phi=\E\left\{ \nabla
V\left(y\right)\cdot\nabla_k\phi_1 +S\left(y\right)\pdr{\phi_1}{z}
\right\}.
\end{equation}
The first term on the right is
\begin{eqnarray}
&&\E\left\{ V\left(y\right)\cdot\nabla_k\phi_1\right\}=
\pdr{}{k_j}\E\left\{V_j(y)\left[\int_0^\infty V_m(y-s k)e^{-\theta
s}ds \pdr{\bar\phi}{k_m}+\int_0^\infty S(y-s k)e^{-\theta
s}ds\pdr{\bar\phi}{z}\right] \right\}\nonumber\\
&&\to \pdr{}{k_j}\left(D_{mj}(k)\pdr{\bar\phi}{k_m}\right)+
\pdr{}{k_j}\left(\int_0^\infty
R_j^{SV}(sk)ds\pdr{\bar\phi}{z}\right)\label{barphi-term1}
\end{eqnarray}
in the limit $\theta\to 0$. The second term in the right side of
(\ref{interm-barphi}) is
\begin{eqnarray}
&&\E\left\{S\left(y\right)\pdr{\phi_1}{z}
\right\}=\pdr{}{z}\E\left\{S(y)\left[\int_0^\infty V_m(y-s
k)e^{-\theta s}ds \pdr{\bar\phi}{k_m}+\int_0^\infty S(y-s
k)e^{-\theta s}ds\pdr{\bar\phi}{z}\right] \right\}\nonumber\\
&&\to -\pdr{}{z}\left(\int_0^\infty
R_m^{VS}(sk)ds\pdr{\bar\phi}{k_m}\right)+ D\frac{\partial^2
\bar\phi}{\partial z^2}=\pdr{}{z}\left(\int_0^\infty
R_m^{SV}(-sk)ds\pdr{\bar\phi}{k_m}\right)+ D\frac{\partial^2
\bar\phi}{\partial z^2}\label{barphi-term2}
\end{eqnarray}
as $\theta\to 0$. Putting together (\ref{barphi-term1}) and
(\ref{barphi-term2}) we obtain
\begin{equation}\label{eq-barphi}
\pdr{\bar\phi}{t}+k\cdot\nabla_x\bar\phi=
\pdr{}{k_j}\left(D_{mj}(k)\pdr{\bar\phi}{k_m}\right)+
\pdr{}{k_j}\left(D_j(k)\pdr{\bar\phi}{z}\right)+\pdr{}{z}\left(
D_j(-k)\pdr{\bar\phi}{k_j}\right)+ D\frac{\partial^2
\bar\phi}{\partial z^2}.
\end{equation}
This is nothing but the Kolmogorov equation for the process with the
generator (\ref{L-comp}).

\section{Proof of Theorem \ref{thm-new}}
\label{sec:proof}

Theorem \ref{thm-new} is a simple corollary of the following
proposition. Let us first introduce some notation. Given a
function $G\in C^{3}_b([0,+\infty)\times\Rm^{2d}_*)$ and $t\ge0$
let us introduce
\[
N_t (G)=G(t,X (t),K (t),Z (t))- G(0,X(0),K(0),Z(0))-
\int\limits_0^t(\partial_\rho  + {\cal L})G(\rho,X (\rho),K
(\rho),Z (\rho))d\rho
\]
with the operator ${\cal L}$ defined in (\ref{L-gen1}). We also denote
by $\Rm^d_*:=\Rm^d\setminus\{0\}$ and $\Rm^{2d}_*:=\Rm^d\times\Rm^d_*$
to avoid the singular point $k=0$. 
We let ${\cal C} :=C([0,+\infty);\Rm^{d}\times\Rm^d_*\times \Rm)$ be
the set of continuous paths of $(X(t),K(t),Z(t))$. For any $u\leq v$
denote by ${\cal M}^{v}_{u}$ the $\sigma$-algebra of subsets of ${\cal
  C}$ generated by $(X(t),K(t),Z(t))$, $t\in[u,v]$. We write ${\cal
  M}^{v}:={\cal M}_{0}^{v}$ and ${\cal M}$ for the $\sigma$ algebra of
Borel subsets of ${\cal C}$. It coincides with the smallest
$\sigma$-algebra that contains all ${\cal M}^{t}$, $t\ge0$. We define
${\cal C}(T,M)$ as the set of paths $\pi\in{\cal C}$ so that both $(2M
)^{-1}\le |K(t)|\le 2M $, and
\begin{eqnarray*}
X(t)-X(u)+\int\limits_u^t  K(s) ds=0 ,\,\mbox{ for all }
\,0\le u<t\le T.\,
\end{eqnarray*}
We denote by ${\mathbb P}_{x,k,z}^\delta$ the probability measure
on ${\cal C}(T,M)$ induced by the trajectories of (\ref{dyn1}) and
by $\E_{x,k,z}^\delta$ the corresponding expectation.
\begin{proposition}
  \label{prop-mart} 
  Suppose that $(x,k)\in{\cal A}(M)=\Rm^d\times\{M^{-1}\le|k|\le M\}$
  for some $M>0$ and that a test function $\zeta\in
  C_b((\Rm^{2d}_*)^{n})$ is non-negative.  Let $\gamma_0\in(0,1/2)$
  and let $0\leq t_1<\cdots< t_n\leq T_*\le t<u\le T$. Assume further
  that $u-t\ge \delta^{\gamma_0}$. Then, there exist constants
  $\gamma_1>0$, $C(T)$ such that for any function $G\in
  C^{3}([T_*,T]\times\Rm^{2d}_*\times\Rm)$ we have
  \begin{equation}
    \label{est-mart} \left|\E_{x,k,z}^\delta\left\{\left[N_u (G)-N_t
   (G)\right] \tilde\zeta\right\}\right| \leq
   C(T)\delta^{\gamma_1}(u-t)\|G\|_{4}
   \Big(\E^{\delta}_{x,k,z}\tilde\zeta + \|\zeta\|_\infty\Big). 
  \end{equation}
  Here $\tilde\zeta(\pi):=\zeta(X(t_1),K(t_1),Z(t_1),\ldots,
  X(t_n),K(t_n),Z(t_n))$, and $\pi=(X(t),K(t),Z(t))$ is any continuous
  path. The choice of the constants $\gamma_1,\,C$ does not depend on
  $(x,k,z)$, $\delta\in(0,1]$, $\zeta$, the times $t_1,\ldots,
  t_n,T_*, T, u,t$, or the test function $G$.
\end{proposition}

{\bf Proof of Theorem \ref{thm-new}}. Theorem \ref{thm-new} is a
simple consequence of Proposition \ref{prop-mart}.  Let
$\phi_0(x,k,z)$ be a test function and let the function
$\bar\phi(t,x,k,z)$ solve the initial value problem
\begin{eqnarray}\label{phi-eq}
&&\pdr{\bar\phi}{t}={\cal L}\bar\phi\\
&&\bar\phi(0,x,k,z)=\phi_0(x,k,z).\nonumber
\end{eqnarray}
We apply Proposition \ref{prop-mart} with
$G(t,x,k,z)=\bar\phi(u-t,x,k,z)$, $t=\delta^\gamma$  and
$u>\delta^\gamma$ with $1/2<\gamma<1$ and take $\tilde\zeta=1$. It
follows from (\ref{est-mart}) that
\begin{equation}\label{appr-mart-sh}
\left|\E_{x,k,z}^\delta\left[\phi_0(X(u),K (u),Z (u))-
\bar\phi(u-\delta^\gamma,X(\delta^\gamma),K (\delta^\gamma), Z
(\delta^\gamma))\right]\right|\le C\|G\|_{4}\delta^{\gamma_1}.
\end{equation}
Using the fact that $\bar\phi$ is a smooth function and
$1/2<\gamma<1$ we conclude that
\begin{equation}\label{appr-mart-sh2}
\left|\E_{x,k,z}^\delta\left[\phi_0(X^\delta(u),K^\delta(u),Z^\delta(u))-
\bar\phi(u,x,k,z)\right]\right|\le
C\|G\|_{4}\delta^{\gamma_1}.
\end{equation}
The conclusion of Theorem \ref{thm-new} now follows. $\Box$

\subsection{The proof of Proposition \ref{prop-mart}}

The proof of Proposition \ref{prop-mart} is technical and uses
the ideas of \cite{BKR-liouv,Kesten-Papanico,KR-Diff,KR-2d}. However, the
present situation is much simpler than in the aforementioned
papers as we already know from Theorem \ref{main-thm} that the
process $(X^\delta,K^\delta)$ converges to a process with the
generator
\[
\tilde{\cal
L}=\pdr{}{k_n}\left(D_{nm}(k)\pdr{}{k_m}\right)-k\cdot\nabla_x.
\]
This means that the law of the process $X^\delta$ will be close to
that of the limit $X(t)$; and in particular $X^\delta(t)$ does not
approach a narrow tube around its past trajectory with a probability
very close to one. This potential return was a major obstacle in the
proofs in \cite{BKR-liouv,Kesten-Papanico,KR-Diff,KR-2d}.

The strategy of the proof is as follows. We will need to deal with
objects of, say, the form
\begin{equation}\label{expr1}
\E\left\{G(X^\delta(s),K^\delta(s),Z^\delta(s))
S\left(\frac{X^\delta(s)}{\delta}\right)
V\left(\frac{X^\delta(s')}{\delta}\right)\right\}
\end{equation}
with $s<s'$ but the difference $s'-s$ small.  Then we will consider a
slight pullback in time
\begin{equation}\label{sigma-def}
\sigma(s)=s-\delta^{1-\gamma},
\end{equation}
with a sufficiently small $\gamma>0$ and a linearization
\begin{equation}\label{linear}
L(\sigma,s)=X^\delta(\sigma)-(s-\sigma)K^\delta(\sigma).
\end{equation}
The characteristic equations (\ref{dyn1}) allow us to estimate the
difference between (\ref{expr1}) and
\begin{equation}\label{expr2}
\E\left\{G(X^\delta(\sigma),K^\delta(\sigma),Z^\delta(\sigma))
S\left(\frac{L(\sigma,s)}{\delta}\right)
V\left(\frac{L(\sigma,s')}{\delta}\right)\right\},
\end{equation}
and show that it is small. However, the latter expectation
approximately splits:
\begin{eqnarray}\label{expr22}
&&\E\left\{G(X^\delta(\sigma),K^\delta(\sigma),Z^\delta(\sigma))
S\left(\frac{L(\sigma,s)}{\delta}\right)
V\left(\frac{L(\sigma,s')}{\delta}\right)\right\}\approx\\
&&\E\left\{G(X^\delta(\sigma),K^\delta(\sigma),Z^\delta(\sigma))\right\}
R^{SV}\left(\frac{L(\sigma,s')-L(\sigma,s)}{\delta}\right).\nonumber
\end{eqnarray}
The reason for the expectation splitting is that the argument of
the function $G$ in (\ref{expr2}) depends only on the potential in
a tube close to the trajectory $X^\delta(t)$ until the the time
$t=\sigma$, while $L(\sigma,s)$ and $L(\sigma,s')$ are at distance
much larger than $\delta$ from this tube with a probability close
to one. Hence, the values of
$G(X^\delta(\sigma),K^\delta(\sigma),Z^\delta(\sigma))$ and, say,
of $V\left(\dfrac{L(\sigma,s')}{\delta}\right)$ are nearly
independent and expectation (\ref{expr2}) splits. This is
formalized by the following mixing lemma.

For any $t\geq0$ we denote by ${\cal F}_t$ the $\sigma$-algebra
generated by $(X^{\delta}(s),K^{\delta}(s),Z^\delta(s))$, $s\leq
t$. Here we suppress, for the sake of abbreviation, writing the
initial data in the notation of the trajectory. We assume that
$X_1,X_2:(\Rm\times\Rm^d\times \Rm^{d^2})^2\rightarrow\Rm$ are
certain continuous functions, $Q$ is a random variable and
$g_1,g_2$ are $\Rm^d$-valued random vectors. We suppose further
that $Q,g_1,g_2$, are ${\cal F}_t$-measurable, while $\tilde
X_1,\tilde X_2$ are random fields of the form
\[
\tilde X_i(x)=X_i \left( F(x), \nabla_x F(x),\nabla_x^2
 F(x)\right)  ,
\]
where, as before we denote for brevity $F=(V,S)$. We also let
\begin{equation}\label{80101}
U(\theta_1,\theta_2):= \E\left[\tilde X_1(\theta_1)\tilde
X_2(\theta_2)\right] ,\quad \theta_1,\theta_2\in\Rm^d
\end{equation}
and recall that $\phi(r)$ is the mixing coefficient defined in
(\ref{mix-coeff}). The following mixing lemma from \cite{BKR-liouv} is
formalizing the expectation splitting.
\begin{lemma}\label{lem-mix}
(i) Assume that $r,t\geq0$ and
\begin{equation} \label{70202}
 \inf\limits_{u\leq t}\left|g_i-\frac{X^{\delta}(u)}{\delta}\right|\geq
\frac{r}{\delta},
\end{equation}
${\mathbb P}$-a.s. on the set $Q\not=0$ for $i=1,2$. Then we have
\begin{equation}
\label{70201} \left|\E\left[\tilde X_1(g_1)\tilde
X_2(g_2)Q\right]-\E\left[ U(g_1,g_2) Q\right]\right| \leq
2\phi\left(\frac{r}{2\delta}\right)\|X_1\|_{L^\infty}
\|X_2\|_{L^\infty}\|Q\|_{L^1(\Omega)}.
\end{equation}
\item(ii) Let $\E X_1(0)=0$ and assume that $g_2$ satisfies
$(\ref{70202})$,
\begin{equation} \label{70202b}
\inf\limits_{u\leq t}\left|g_1
-\frac{X^{\delta}(u)}{\delta}\right|\geq \frac{r+r_1}{\delta}
\end{equation}
and $|g_1 -g_2 |\geq r_1\delta^{-1}$ for some $r_1\geq0$, ${\mathbb
P}$-a.s. on the event $Q\not=0$.  Then, we have
\begin{equation}
\label{70201b} \left|\E\left[\tilde X_1(g_1)\tilde
X_2(g_2)\,Q\right]-\E\left[ U(g_1,g_2) Q\right]\right| \leq
C\phi^{1/2}\left(\frac{r}{2\delta}\right)
\phi^{1/2}\left(\frac{r_1}{2\delta}\right)
\|X_1\|_{L^\infty}\|X_2\|_{L^\infty}\|Q\|_{L^1(\Omega)}
\end{equation}
for some absolute constant $C>0$.
Here the function $U$ is given by $(\ref{80101})$.
\end{lemma}

We proceed now with the proof of Proposition \ref{prop-mart}. Let
$G(k,z)$ be a sufficiently smooth function. We will establish the
approximate martingale property (\ref{est-mart}) for $G$. It
suffices to consider functions of the form $G(z,k)=g(z)r(k)$. The
characteristic equations
\[
Z^\delta(t)=z+\frac{1}{\sqrt{\delta}}\int_0^t
S\left(\frac{X^\delta(s)}{\delta}\right)ds,
\]
and
\[
K_j^\delta(t)=k+\frac{1}{\sqrt{\delta}}\int_0^t
V_j\left(\frac{X^\delta(s)}{\delta}\right)ds,
\]
imply that we have
\begin{eqnarray}\label{d1}
&&G(K^\delta(u),Z^\delta(u))-G(K^\delta(t),Z^\delta(t))\\
&&= \frac{1}{\sqrt{\delta}}\int_t^u
\left[{g'(Z^\delta(s))}S\left(\frac{X^\delta(s)}{\delta}\right)r(K^\delta(s))+
g(Z^\delta(s))r_j(K^\delta(s))V_j\left(\frac{X^\delta(s)}{\delta}\right)
\right]ds.
\nonumber
\end{eqnarray}
Here and below, we use the notation $V_j$ for $V_{x_j}$ and similarly
$S_k$ for $S_{x_k}$ to simplify.  In order to be able to make a
backward step $\sigma$ as in (\ref{sigma-def}) we split the above
integral as
\begin{eqnarray}\label{c2}
&&G(K^\delta(u),Z^\delta(u))-G(K^\delta(t),Z^\delta(t))= A+B
\end{eqnarray}
with
\[
A=\frac{1}{\sqrt{\delta}}\int_t^{t+\delta^{1-\gamma}} \left[
g'(Z^\delta(s))r(K^\delta(s))
S\left(\frac{X^\delta(s)}{\delta}\right)+
g(Z^\delta(s))r_j(K^\delta(s))V_j\left(\frac{X^\delta(s)}{\delta}\right)\right]
ds
\]
and
\[
B=\frac{1}{\sqrt{\delta}}\int_{t+\delta^{1-\gamma}}^u
 \left[g'(Z^\delta(s))r(K^\delta(s))
S\left(\frac{X^\delta(s)}{\delta}\right)+
g(Z^\delta(s))r_j(K^\delta(s))V_j\left(\frac{X^\delta(s)}{\delta}\right)\right]ds=
B_1+B_2.
\]
The first term is small:
\[
|A|\le C\delta^{1/2-\gamma}\|G\|_{2}
\]
provided that $\gamma<1/2$. The term $B$ will be analyzed with the
aforementioned ideas of a backward step and linearization, and
using the mixing lemma. The terms $B_1$ and $B_2$ may be written
as
\begin{eqnarray}\nonumber
&&B_1=\frac{1}{\sqrt{\delta}}\int_{t+\delta^{1-\gamma}}^u
\left\{g'(Z^\delta(\sigma(s)))r(K^\delta(\sigma(s))+
[g'(Z^\delta(s))r(K^\delta(s))-g'(Z^\delta(\sigma(s)))r(K^\delta(\sigma(s))]\right\}
\\
&&~~~~~~~~~~\times
S\left(\frac{X^\delta(s)}{\delta}\right)ds=I+II\label{d21}
\end{eqnarray}
and
\begin{eqnarray}\nonumber
&&B_2=\frac{1}{\sqrt{\delta}}\int_{t+\delta^{1-\gamma}}^u
\left\{g(Z^\delta(\sigma(s)))r_j(K^\delta(\sigma(s))+
[g(Z^\delta(s))r_j(K^\delta(s))-g(Z^\delta(\sigma(s)))r_j(K^\delta(\sigma(s))]\right\}
\\
&&~~~~~~~~~~\times
V_j\left(\frac{X^\delta(s)}{\delta}\right)ds=III+IV\label{d22}
\end{eqnarray}
 with
\begin{equation}\label{d3}
I=\frac{1}{\sqrt{\delta}}\int_{t+\delta^{1-\gamma}}^u
g'(Z^\delta(\sigma(s)))r(K^\delta(\sigma(s)))
S\left(\frac{X^\delta(s)}{\delta}\right)ds
\end{equation}
and
\begin{eqnarray}
&& II=\frac{1}{\sqrt{\delta}}\int_{t+\delta^{1-\gamma}}^u
[g'(Z^\delta(s))r(K^\delta(s))-g'(Z^\delta(\sigma(s)))r(K^\delta(\sigma(s))]
S\left(\frac{X^\delta(s)}{\delta}\right)ds \label{d4}
\end{eqnarray}
while
\begin{eqnarray}\label{d5}
III=\frac{1}{\sqrt{\delta}}\int_{t+\delta^{1-\gamma}}^u
g(Z^\delta(\sigma(s)))r_j(K^\delta(\sigma(s))
V_j\left(\frac{X^\delta(s)}{\delta}\right)ds
\end{eqnarray}
and
\begin{eqnarray}
&&IV=\frac{1}{\sqrt{\delta}}\int_{t+\delta^{1-\gamma}}^u
[g(Z^\delta(s))r_j(K^\delta(s))-g(Z^\delta(\sigma(s)))r_j(K^\delta(\sigma(s))]
V_j\left(\frac{X^\delta(s)}{\delta}\right)ds.\label{d6}
\end{eqnarray}
We summarize the contributions of each of the terms above in the
following lemma.
\begin{lemma}\label{lem-I-IV}
There exists a constant $\alpha>0$ so that
the terms $I$, $II$, $III$ and $IV$ satisfy the following
estimates:
\begin{eqnarray}\label{I-est-lem}
&&\left|\E\left\{\left[I -\int_{t}^u\int_0^\infty \theta \Delta
R^{SV}\left(\theta K^\delta(s) \right)d\theta g'(Z^\delta(
s))r(K^\delta(s))ds \right]\tilde\zeta\right\}\right|\\
&&\le C\delta^{\alpha}\|G\|_{4}\left[
\|\zeta\|_\infty+\E\left\{\tilde\zeta\right\}\right](u-t),\nonumber\\
\label{II-est-lem} &&\left|\E\left[\left\{II-\int_t^u\int_0^\infty
g''(Z^\delta(s))r(K^\delta(s))
 R^{SS}(\theta K^\delta(s))d\theta
ds\right.\right.\right.\\
&&\left.\left.\left.-\int_t^u\int_0^\infty
g'(Z^\delta(s))r_j(K^\delta(s)) R_j^{SV}(\theta
K^\delta(s))d\theta ds\right\} \tilde\zeta\right]\right|\le
C\delta^{\alpha}\|G\|_3\left[\|\zeta\|_\infty+
E\left\{\tilde\zeta\right\}\right](u-t),\nonumber\\
\label{III-est-lem} &&\left|\E\left\{III \tilde\zeta \right\}+
 \int_{t }^u\int_0^\infty
\theta g(Z^\delta(s))r_j(K^\delta(s)) \Delta
R_{j}^{VV}\left(\theta
K^\delta(s)\right) d\theta ds\right|\\
 &&\le C\delta^\alpha\|G\|_1\left(\|\zeta\|_\infty+\E[\tilde\zeta]\right)(u-t)
 \nonumber\\
 \label{IV-est-lem}
&&\left|\E\{IV\tilde\zeta\}-
\int\limits_{t}^u\int\limits_{0}^\infty
g'(Z^\delta(s))r_j(K^\delta(s)) R_j^{SV}\left(-\theta K^\delta(s)
\right)d\theta ds \right.  \\
&&\left. + \int\limits_{t}^u\int\limits_{0}^\infty
g(Z^\delta(s))r_{jm}(K^\delta(s))R_{mj}^{VV}\left(\theta
K^\delta(s) \right) d\theta ds\tilde\zeta\right| \le
C\delta^\alpha\|G\|_2\left(\|\tilde\zeta\|_\infty+E\{\tilde\zeta\}\right)(u-t).
\nonumber
\end{eqnarray}
\end{lemma}
It remains now only to prove Lemma \ref{lem-I-IV} as the four
individual contributions above combine to the operator ${\cal L}$
in (\ref{L-gen1}).

\subsection{The estimate for $I$}

We first recall the linear approximation (\ref{linear}) and define
the interpolation
\[
R(v,\sigma,s)=(1-v)L(\sigma,s)+vX^\delta(s).
\]
This allows us to write a linear approximation for $S$ as
\begin{eqnarray*}
&&S\left(\frac{X^\delta(s)}{\delta}\right)=
S\left(\frac{R(1,\sigma,s)}{\delta}\right)=S\left(\frac{R(0,\sigma,s)}{\delta}\right)+
\frac{1}{\delta}\int_0^1S_{i}\left(\frac{R(v,\sigma,s)}{\delta}\right)(X_i^\delta(s)-
L_i(\sigma,s))dv\\
&& ~~~~~~~~~~~~~~~~=S\left(\frac{L(\sigma,s)}{\delta}\right)+
\frac{1}{\delta}\int_0^1S_{i}\left(\frac{R(v,\sigma,s)}{\delta}\right)(X_i^\delta(s)-
L_i(\sigma,s))dv.
\end{eqnarray*}
Now we split $I$ as
\begin{equation}\label{IJ}
I=J_1+J_2
\end{equation}
according to the above, with
\begin{eqnarray}\label{J1}
J_1=\frac{1}{\sqrt{\delta}}\int_{t+\delta^{1-\gamma}}^u
g'(Z^\delta(\sigma(s)))r(K^\delta(\sigma(s))S\left(\frac{L(\sigma,s)}{\delta}\right)ds
\end{eqnarray}
and
\begin{eqnarray}\label{J2}
J_2=\frac{1}{ {\delta}^{3/2}}\int_{t+\delta^{1-\gamma}}^u\int_0^1
g'(Z^\delta(\sigma(s)))r(K^\delta(\sigma(s)))
S_{i}\left(\frac{R(v,\sigma,s)}{\delta}\right)(X_i^\delta(s)-
L_i(\sigma,s))dvds.
\end{eqnarray}
The term $J_1$ is ready for an application of the mixing lemma: the
arguments of the function $G$ (that is, of $g$ and $r$), and of the
field $S$ are separated by a distance of the order
$O(\delta^{1-\gamma})$ that is much larger than $\delta$, with a
probability close to one. In order to make this statement precise we
introduce a stopping time $\tau_\delta$ that ensures that until
$\tau_\delta$ the trajectory $X^\delta(t)$ "goes forward" and does not
come back to its past.

Let $0<\eps_1<\eps_2<1/2$, $\eps_3\in(0,1/2-\eps_2)$,
$\eps_4\in(1/2,1-\eps_1-\eps_2)$ be small positive constants and
set
\begin{equation}\label{102302}
N=[\delta^{-\eps_1}],\quad p=[\delta^{-\eps_2}],\quad
q=p\,[\delta^{-\eps_3}],\quad N_1=Np\,[\delta^{-\eps_4}].
\end{equation}
The requirement is that $\eps_i$, $i\in\{1,2,3\}$ should be
sufficiently small and $\eps_4$ is bigger than $1/2$, less than one
and can be made as close to one as we would need it. It is important
that $\eps_1<\eps_2$ so that $N\ll p$ when $\delta\ll 1$. We introduce
the following $({\cal M}^{t})_{t\geq0}$ stopping times. Let
$t^{(p)}_k:=kp^{-1}$ be a mesh of times, and $\pi\in {\cal C}$ be a
path.  We define the ``violent turn'' stopping time
\begin{eqnarray}\label{Sdelta}
&&V_\delta(\pi):=\inf\left[\,t\geq0:\vphantom{\int_0^1}\mbox{ for
some }k\geq 0\mbox{ we have
}t\in\left[t_k^{(p)},t_{k+1}^{(p)}\right)\right. \mbox{ and }\\
&&\left. ~~~~~~~ \hat{K}(t_{k-1}^{(p)})\cdot \hat K(t)\le
1-\frac{1}{N},\mbox{ or }\,
\hat{K}\left(t_{k}^{(p)}-\frac{1}{N_1}\right)\cdot \hat
K(t)\le1-\frac{1}{N}\,\right],\nonumber
\end{eqnarray}
where by convention we set $\hat K(-1/p):=\hat K(0)$.  Note that with
the above choice of $\eps_4$ we have
$\hat{K}\left(t_{k}^{(p)}-1/N_1\right)\cdot \hat
K(t_{k}^{(p)})>1-1/N$, provided that $\delta\in(0,\delta_0]$ and
$\delta_0$ is sufficiently small.  The stopping time $V_\delta$ is
triggered when the trajectory performs a sudden turn; this is
undesirable as the trajectory may then return back to the region it
has already visited and create correlations with the past.

For each $t\ge 0$, we denote by
 $ \mathfrak X_t(\pi):=\mathop{\bigcup\limits_{0\le s\leq
t}}X\left(s;\pi\right)$ the trace of the spatial component of the
path $\pi$ up to time $t$, and by $\mathfrak
X_{t}(q;\pi):=[\bx:\mbox{dist }(\bx, \mathfrak X_{t}(\pi))\le
1/q]$ a tubular region around the path. We introduce the stopping
time
\begin{equation}\label{Udelta}
U_\delta(\pi):=\inf\left[\,t\ge0:\,\exists\,k\ge1\, \mbox{ and
}t\in[t_k^{(p)},t_{k+1}^{(p)})\mbox{ for which }X(t)\in \mathfrak
X_{t_{k-1}^{(p)}}(q)\,\right].
\end{equation}
It is associated with the return of the $X$ component of the
trajectory to the tube around its past; this is again an undesirable
way to create correlations with the past. Finally, we set the stopping
time
\begin{equation}
\label{W:delta}
 \tau_\delta(\pi):=V_\delta(\pi)\wedge U_\delta(\pi).
\end{equation}
\begin{lemma}\label{lem-tau}\cite{KR-Diff}
The probability of the event $[\,\tau_\delta< T\,]$ for a fixed
$T>0$ goes to zero, as $\delta\to 0$: there exists $\alpha_0>0$ so
that
\begin{equation}\label{tau-prob}
{\mathbb P}\left\{[\,\tau_\delta< T\,]\right\}\le
C(T)\delta^{\alpha_0}.
\end{equation}
\end{lemma}

We apply part (i) of Lemma \ref{lem-mix} to $E\{J_1\tilde\zeta\}$
with 
\begin{displaymath}
  \tilde X_1(x)= S(x),\quad \tilde X_2=1,\quad
  {Q}=g'\big(Z^\delta(\sigma)\big)r\big(K^\delta(\sigma(s))\big){\bf
   1}[\tau_\delta>T]\tilde\zeta, \quad
  g_1= \dfrac{L (\sigma,s)}{\delta}.
\end{displaymath}
Note that $g_1$ and $Q$ are both
${\cal F}_\sigma$ measurable. It follows from the definition of
the stopping time $\tau_\delta$ that when ${Q}\ne 0$ then the
linearization also stays away from the past trajectory: for all
$0\le\rho\le\sigma$ we have
\begin{equation}\label{L-X-far}
\left|L(\sigma,s)-X^\delta(\rho)\right|\ge C\delta^{1-\gamma}
\end{equation}
and hence
\[
\inf_{0\le\rho\le\sigma(s)}\left|g_1-\frac{X^\delta(\rho)}{\delta}\right|\ge
\frac{r}{\delta}
\]
with $r=C\delta^{1-\gamma}$. We also decompose $J_1$ according to
whether the stopping time has occurred or not before time $T$:
\begin{eqnarray*}
&&J_1=\frac{1}{\sqrt{\delta}}\int_{t+\delta^{1-\gamma}}^u
g'(Z^\delta(\sigma(s)))r(K^\delta(\sigma(s))) {\bf
1}[\tau_\delta>T]S\left(\frac{L(\sigma,s)}{\delta}\right)ds\\
&&~~~~+ \frac{1}{\sqrt{\delta}}\int_{t+\delta^{1-\gamma}}^u
g'(Z^\delta(\sigma(s)))r(K^\delta(\sigma(s))) (1-{\bf
1}[\tau_\delta>T])S\left(\frac{L(\sigma,s)}{\delta}\right)ds
=J_{11}+J_{12}.
\end{eqnarray*}
However, (\ref{tau-prob}) implies that
\begin{equation}\label{dJ12}
\E\left\{\left|J_{12}\tilde\zeta\right|\right\}\le
C\delta^{\alpha_0}\|G\|_{2}\|\zeta\|_\infty(u-t)
\end{equation}
so we have to deal only with $J_{11}$. Using the mixing lemma as
above, with the point separation as in (\ref{L-X-far}), and the fact
that $\E[S(x)]=0$ (whence $U=0$ in \eqref{70201}) we estimate
\begin{eqnarray}\label{dc5}
\left|\E\left(J_{11}\tilde\zeta\right)\right|\le\frac{C}{\sqrt{\delta}}
\phi\left(C\delta^{-\gamma}\right)(u-t)\|g\|_1\|r\|_0\E[\tilde\zeta]\le
C\delta\|G\|_{2}\E[\tilde\zeta](u-t).
\end{eqnarray}
Here $\phi(r)$ is the mixing coefficient that decays faster than any
power of $r$; see (\ref{DR}). We conclude that
\begin{equation}\label{d-Jest}
\left|\E\left(J_{1}\tilde\zeta\right)\right|\le
C\delta\|G\|_{2}\E[\tilde\zeta](u-t)\|
+C\delta^{\alpha_0}\|G\|_{2}\zeta\|_\infty(u-t)
\end{equation}
so the term $J_1$ produces only a small contribution.

Now we estimate the second term $J_2$ in (\ref{IJ}); it is given
explicitly by (\ref{J2}). We split it further by using the next order
expansion
\[
S_{i}\left(\frac{R(v,\sigma,s)}{\delta}\right)=
S_{i}\left(\frac{L(\sigma,s)}{\delta}\right)+
\frac{1}{\delta}\int_0^v
S_{ij}\left(\frac{R(\theta,\sigma,s)}{\delta}\right)
(X_j^\delta(s)-L_j(\sigma,s))d\theta.
\]
This leads to the corresponding expression $J_2=J_{21}+J_{22}$
with
\begin{eqnarray*}
J_{21}=\frac{1}{ {\delta}^{3/2}}\int_{t+\delta^{1-\gamma}}^u
g'(Z^\delta(\sigma(s)))r(K^\delta(\sigma(s)))
S_{i}\left(\frac{L(\sigma,s)}{\delta}\right)(X_i^\delta(s)-
L_i(\sigma,s))ds
\end{eqnarray*}
and
\begin{eqnarray*}
&&J_{22}=\frac{1}{
{\delta}^{5/2}}\int\limits_{t+\delta^{1-\gamma}}^u
\int\limits_0^1\int\limits_0^v
g'(Z^\delta(\sigma(s)))r(K^\delta(\sigma(s)))
S_{ij}\left(\frac{R(\theta,\sigma,s)}{\delta}\right)\\
&&~~~~~~~~~\times(X_i^\delta(s)-
L_i(\sigma,s))(X_j^\delta(s)-L_j(\sigma,s))d\theta dvds.
\end{eqnarray*}
Note that the characteristic equations and the definition
(\ref{sigma-def}) of the time $\sigma(s)$ imply that
\begin{equation} \label{L-X-close}
|L(\sigma,s)-X^\delta(s)|\le
C\delta^{2-2\gamma-1/2}=C\delta^{3/2-2\gamma}.
\end{equation}
It follows that
\begin{equation}\label{J22-small}
\left|\E\left\{J_{22}\tilde\zeta\right\}\right|\le C\|G\|_{2}\E
\{\tilde\zeta \}\delta^{-5/2}\delta^{3-4\gamma}(u-t)\le
C\delta^{1/2-4\gamma}\|G\|_{2}\E \{\tilde\zeta \}(u-t),
\end{equation}
which is small if $\gamma<1/8$. This is an important characteristic
feature of the weak coupling limit; after several linearizations the
remainder becomes deterministically small while the linearized terms
may be controlled with the mixing lemma.

Next, we look at $J_{21}$ that is the only potentially surviving
(not small) in the limit $\delta\to 0$ contribution to $I$: to do
so we write, using the evolution equation for $X^\delta$ and a
further linearization for the function $V_i$:
\begin{eqnarray}
&&\!\!\!\!\!\!\!\!\!\!\!\!\!\!\!\!\!\!\!\!\!\!\!\!\!\!
X_i^\delta(s)-L_i(\sigma,s)=-\int_\sigma^s[K_i^\delta(v)-K_i^\delta(\sigma)]dv=
-\frac{1}{\sqrt{\delta}}\int_\sigma^s\int_\sigma^{v}
V_i\left(\frac{X^\delta(\rho)}{\delta}\right)d\rho dv\nonumber\\
&& = -\frac{1}{\sqrt{\delta}}\int_\sigma^s(s-\rho)
V_i\left(\frac{X^\delta(\rho)}{\delta}\right)d\rho=-
\frac{1}{\sqrt{\delta}}\int_\sigma^s(s-\rho)
V_i\left(\frac{L(\sigma,\rho)}{\delta}\right)d\rho\nonumber\\
&&-\frac{1}{\sqrt{\delta}}\int_\sigma^s(s-\rho)
\left[V_i\left(\frac{X^\delta(\rho)}{\delta}\right)-
V_i\left(\frac{L(\sigma,\rho)}{\delta}\right)\right]d\rho=
-\frac{1}{\sqrt{\delta}}\int_\sigma^s(s-\rho)
V_i\left(\frac{L(\sigma,\rho)}{\delta}\right)d\rho\nonumber\\
&&-\frac{1}{\delta^{3/2}}\int_\sigma^s\int_0^1(s-\rho)
V_{im}\left(\frac{R(v,\sigma,\rho)}{\delta}\right)
\left(X_m^\delta(\rho)-L_m(\sigma,\rho)\right)dvd\rho.\label{L-X-approx}
\end{eqnarray}
This produces a further split
\[
J_{21}=J_{21}^1+J_{21}^2
\]
with
\begin{eqnarray*}
J_{21}^1=-\frac{1}{
{\delta}^{2}}\int_{t+\delta^{1-\gamma}}^u\int_\sigma^s(s-\rho)
g'(Z^\delta(\sigma(s)))r(K^\delta(\sigma(s)))
S_{i}\left(\frac{L(\sigma,s)}{\delta}\right)
V_i\left(\frac{L(\sigma,\rho)}{\delta}\right)d\rho ds
\end{eqnarray*}
and
\begin{eqnarray*}
&&J_{21}^2=-\frac{1}{
{\delta}^{3}}\int\limits_{t+\delta^{1-\gamma}}^u
\int\limits_\sigma^s\int\limits_0^1 (s-\rho)
g'(Z^\delta(\sigma(s)))r(K^\delta(\sigma(s)))
S_{i}\left(\frac{L(\sigma,s)}{\delta}\right)
V_{im}\left(\frac{R(v,\sigma,\rho)}{\delta}\right)\\
&&~~~~~~~~\times\left(X_m^\delta(\rho)-L_m(\sigma,\rho)\right)dvd\rho
ds.
\end{eqnarray*}
Note that we have linearized enough to  achieve a deterministic
estimate
\begin{equation}\label{J212-small}
|J_{21}^2|\le
C\delta^{-3}\delta^{1-\gamma}\delta^{1-\gamma}\delta^{3/2-\gamma}\|G\|_1(u-t)=
C\delta^{1/2-3\gamma}\|G\|_{1}(u-t)
\end{equation}
which is small if $\gamma<1/6$. Hence, $J_{21}^1$ remains the only
potentially contributing term in $I$: it is analyzed with the help
of the mixing lemma. We take
\[
g_1=\frac{L(\sigma,s)}{\delta},~~g_2=\frac{L(\sigma,\rho)}{\delta},~~
X_1(x)=S_{i}(x),~~X_2(x)=V_i(x),~~r=\rho-\sigma,~~r_1=s-\rho,
\]
and
\[
{Q}=g'(Z^\delta(\sigma))r(K^\delta(\sigma(s))){\bf
1}[\tau_\delta>T]\tilde\zeta.
\]
Let us denote
\[
R^{SV}_{ij}(x)=-\E\left\{S_i(y)V_j(x+y)\right\}.
\]
Observe that
\[
R^{SV}_{ij}(x)=-\E\left\{S_i(y)V_j(x+y)\right\}=-\pdr{}{x_j}
\E\left\{S_i(y)V(x+y)\right\}=
\pdr{}{x_j}\E\left\{S(y)V_i(x+y)\right\}=
\frac{\partial^2R^{SV}(x)}{\partial x_i\partial x_j}.
\]
The mixing lemma implies that
\begin{eqnarray*}
&&\!\!\!\!\left|\E\left[J_{21}^1\tilde\zeta\right]+\frac{1}{
{\delta}^{2}}\int\limits_{t+\delta^{1-\gamma}}^u\int\limits_\sigma^s(s-\rho)
\E\left\{g'(Z^\delta(\sigma(s)))r(K^\delta(\sigma(s)))
\left(-R_{ii}^{SV}\left(\frac{L(\sigma,\rho)-L(\sigma,s)}{\delta}
\right)\right)\right\}d\rho ds\right|\\
&&\le\frac{C\|G\|_{2}\E(\tilde\zeta)}{\delta^2}
\int\limits_{t+\delta^{1-\gamma}}^u\int\limits_\sigma^s(s-\rho)
\phi^{1/2}\left(\frac{\rho-\sigma}{\delta}\right)
\phi^{1/2}\left(\frac{s-\rho}{\delta}\right)d\rho ds\le
C\delta\|G\|_{2}\E(\tilde\zeta)(u-t)
\end{eqnarray*}
since the mixing coefficient $\phi$ is rapidly decaying. It
follows that
\begin{eqnarray*}
&&\!\!\left|\E\left[J_{21}^1\tilde\zeta\right]-\frac{1}{
{\delta}^{2}}\int_{t+\delta^{1-\gamma}}^u\int_\sigma^s(s-\rho)
\E\left\{g'(Z^\delta(\sigma(s)))r(K^\delta(\sigma(s)))\left(\Delta
R^{SV}\left(\frac{(s-\rho)K^\delta(\sigma)}{\delta}
\right)\right)\right\}d\rho ds\right|\\
&& \le C\delta\|G\|_{1,1}\E(\tilde\zeta)(u-t)
\end{eqnarray*}
The integral above from $\sigma$ to $s$ may be massaged as
\begin{eqnarray*}
&& \frac{1}{ {\delta}^{2}} \int_\sigma^s(s-\rho)
 \left(\Delta R^{SV}\left(\frac{(s-\rho)K^\delta(\sigma)}{\delta}
\right)\right) d\rho = \int_0^{\delta^{-\gamma}}\theta \Delta
R^{SV}\left(\theta K^\delta(\sigma) \right)d\theta.
\end{eqnarray*}
Observe that
\[
|K^\delta(s)-K^\delta(\sigma)|\le C\delta^{1/2-\gamma}
\]
and the function $R^{SV}$ is smooth and rapidly decaying. This allows
us to replace $\sigma$ in the argument of $R^{SV}$ by $s$.  The same
can be done with the functions $g$ and $r$; we conclude that
\begin{equation}\label{d-J211}
\left|\E\left\{\left[J_{21}^1 -\int_{t}^u\int_0^\infty \theta
\Delta R^{SV}\left(\theta K^\delta(s) \right)d\theta g'(Z^\delta(
s))r(K^\delta(s))ds \right]\tilde\zeta\right\}\right|\le
C\delta^{\gamma_1}\|G\|_{4}\E\left\{\tilde\zeta\right\}(u-t)
\end{equation}
with some $\gamma_1>0$. Therefore, putting all the work of this
section together, see (\ref{d-Jest}), (\ref{J22-small}) and
(\ref{J212-small}), we also have
\begin{equation}\label{d-I-est}
\left|\E\left\{\left[I -\int_{t}^u\int_0^\infty \!\!\theta \Delta
R^{SV}\!\left(\theta K^\delta(s) \right)d\theta g'(Z^\delta(
s))r(K^\delta(s))ds \right]\tilde\zeta\right\}\right|\le
C\delta^{\gamma_1}\|G\|_{4}\left[
\|\zeta\|_\infty+\E\left\{\tilde\zeta\right\}\right](u-t).
\end{equation}
This proves the estimate (\ref{I-est-lem}) in Lemma
\ref{lem-I-IV}.

\subsection{Estimate for $II$}

We now look at the term $II$ given by (\ref{d4}) and split it
further as
\begin{eqnarray*}
&&II=\frac{1}{\sqrt{\delta}}\int_{t+\delta^{1-\gamma}}^u
[g'(Z^\delta(s))r(K^\delta(s))-g'(Z^\delta(\sigma(s)))r(K^\delta(\sigma(s))]
S\left(\frac{X^\delta(s)}{\delta}\right)ds \\
&&~~~~=\frac{1}{\delta}
\int\limits_{t+\delta^{1-\gamma}}^u\int\limits_{\sigma(s)}^s
g''(Z^\delta(\rho))r(K^\delta(\rho))S\left(\frac{X^\delta(\rho)}{\delta}\right)
S\left(\frac{X^\delta(s)}{\delta}\right){d\rho ds}\label{d41}\\
&&~~~~+\frac{1}{\delta}
\int\limits_{t+\delta^{1-\gamma}}^u\int\limits_{\sigma(s)}^s
g'(Z^\delta(\rho))r_j(K^\delta(\rho))
V_j\left(\frac{X^\delta(\rho)}{\delta}\right)
S\left(\frac{X^\delta(s)}{\delta}\right){d\rho ds},
\end{eqnarray*}
The estimation of $II$ is very similar to that of $I$ both in
spirit and in mechanics but is even somewhat simpler since as all
we have to justify is the replacement of the arguments of $S$ and
$V_j$ by the corresponding linear approximation from the time
$\sigma(s)$. This is done as in the previous section with the help
of the mixing lemma and linearization and leads to
\begin{equation}\label{II-II'}
\left|\E\left\{II\tilde\zeta\right\}-\E\left\{II'\tilde\zeta\right\}\right|\le
C\delta^{\gamma_2}\|G\|_{3}\left[\|\zeta\|_\infty+\E\{\tilde\zeta\}\right](u-t)
\end{equation}
with
\begin{eqnarray*}
&&II'= \frac{1}{\delta}
\int\limits_{t+\delta^{1-\gamma}}^u\int\limits_{\sigma(s)}^s
g''(Z^\delta(\sigma))r(K^\delta(\sigma))S\left(\frac{L(\sigma,\rho)}{\delta}\right)
S\left(\frac{L(\sigma,s)}{\delta}\right){d\rho ds}\label{d411}\\
&&~~~~+\frac{1}{\delta}
\int\limits_{t+\delta^{1-\gamma}}^u\int\limits_{\sigma(s)}^s
g'(Z^\delta(\sigma))r_j(K^\delta(\sigma))
V_j\left(\frac{L(\sigma,\rho)}{\delta}\right)
S\left(\frac{L(\sigma,s)}{\delta}\right){d\rho ds}=II_1+II_2.
\end{eqnarray*}
The mixing lemma and rapid decay of the mixing coefficient imply
that
\[
\left|\E\left[\left\{II_1-\frac{1}{\delta}
\int\limits_{t+\delta^{1-\gamma}}^u\int\limits_{\sigma }^s\!
g''(Z^\delta(\sigma))r(K^\delta(\sigma))
R^{SS}\!\left(\frac{L(\sigma,s)-L(\sigma,\rho)}{\delta}\right)
d\rho ds\right\}\tilde\zeta\right]\right|\le C\delta
\|G\|_3E\{\tilde\zeta\}(u-t)
\]
The inner integral may be re-written as
\begin{eqnarray*}
&&\!\!\!\!\!\int\limits_{\sigma}^s\! g''(Z^\delta(\sigma))r(K^\delta(\sigma))
R^{SS}\!\left(\frac{L(\sigma,s)-L(\sigma,\rho)}{\delta}\right)
d\rho
=\int\limits_{\sigma}^s \!g''(Z^\delta(\sigma))r(K^\delta(\sigma))
R^{SS}\left(\frac{(\rho-s)K^\delta(\sigma)}{\delta}\right)
d\rho\\
&&=
\int_0^{\delta^{-\gamma}}g''(Z^\delta(\sigma))r(K^\delta(\sigma))
R^{SS}(\theta K^\delta(\sigma))d\theta
\end{eqnarray*}
Using the rapid decay of $R^{SS}$, smoothness of $G$ and closeness
of $s$ and $\sigma$ we obtain
\begin{equation}\label{d-II1-est}
\left|\E\left[\left\{II_1-\int_t^u\int_0^\infty
g''(Z^\delta(s))r(K^\delta(s))
 R^{SS}(\theta K^\delta(s))d\theta
ds\right\}\tilde\zeta\right]\right|\le
C\delta^{\gamma_2}\|G\|_3E\left\{\tilde\zeta\right\}(u-t).
\end{equation}

Similarly, we have for $II_2$:
\[
\left|\E\left[\left\{II_2+\frac{1}{\delta}
\int\limits_{t+\delta^{1-\gamma}}^u\int\limits_{\sigma }^s
g'(Z^\delta(\sigma))r_j(K^\delta(\sigma))
R_j^{VS}\!\!\left(\frac{L(\sigma,s)-L(\sigma,\rho)}{\delta}\right)
\!d\rho ds\right\}\tilde\zeta\right]\right|\le C\delta
\|G\|_3E\{\tilde\zeta\}(u-t)
\]
The inner integral may be re-written as
\begin{eqnarray*}
&&\!\!\!\!\!\!\int\limits_{\sigma}^s\! g'(Z^\delta(\sigma))r_j(K^\delta(\sigma))
R_j^{VS}\!\left(\frac{L(\sigma,s)-L(\sigma,\rho)}{\delta}\right)
d\rho
\!=\!\int\limits_{\sigma}^s\! g'(Z^\delta(\sigma))r_j(K^\delta(\sigma))
R_j^{VS}\!\left(\frac{(\rho-s)K^\delta(\sigma)}{\delta}\right)
d\rho\\
&&=
\int_0^{\delta^{-\gamma}}g'(Z^\delta(\sigma))r_j(K^\delta(\sigma))
R_j^{VS}(-\theta
K^\delta(\sigma))d\theta=-\int_0^{\delta^{-\gamma}}
g'(Z^\delta(\sigma))r_j(K^\delta(\sigma)) R_j^{SV}(\theta
K^\delta(\sigma))d\theta
\end{eqnarray*}
Using the rapid decay of $R^{SS}$, smoothness of $G$ and closeness
of $s$ and $\sigma$ we obtain
\begin{equation}\label{d-II2-est}
\left|\E\left[\left\{II_2-\int_t^u\int_0^\infty
g'(Z^\delta(s))r_j(K^\delta(s)) R_j^{SV}(\theta
K^\delta(s))d\theta ds\right\}\tilde\zeta\right]\right|\le
C\delta^{\gamma_2}\|G\|_3\left[
\|\zeta\|_\infty+E\left\{\tilde\zeta\right\}\right](u-t).
\end{equation}
Putting (\ref{II-II'}), (\ref{d-II1-est}) and (\ref{d-II2-est})
together we obtain (\ref{II-est-lem}).

\subsection{Estimate of $III$}

We look at the third term in (\ref{d22}) given by (\ref{d5})
\begin{eqnarray}\label{d51}
III=\frac{1}{\sqrt{\delta}}\int_{t+\delta^{1-\gamma}}^u
g(Z^\delta(\sigma(s)))r_j(K^\delta(\sigma(s))
V_j\left(\frac{X^\delta(s)}{\delta}\right)ds.
\end{eqnarray}
This term is similar to $I$ in that before using the mixing lemma
we have to expand a little bit:
\begin{eqnarray*}
&&V_j\left(\frac{X^\delta(s)}{\delta}\right)=
V_j\left(\frac{R(1,\sigma,s)}{\delta}\right)=
V_j\left(\frac{R(0,\sigma,s)}{\delta}\right)+
\int_0^1\frac{d}{dv}V_j\left(\frac{R(v,\sigma,s)}{\delta}\right)dv\\
&&~~~~~~~~~~~~~~~~= V_j\left(\frac{L(\sigma,s)}{\delta}\right)+
\frac{1}{\delta}\int_0^1V_{jk}\left(\frac{R(v,\sigma,s)}{\delta}\right)
(X_k^\delta(s)-L_k(\sigma,s))dv.
\end{eqnarray*}
Accordingly we split $III$ as $III=III_1+III_2$ with
\begin{eqnarray}\label{d61}
III_1=\frac{1}{\sqrt{\delta}}\int_{t+\delta^{1-\gamma}}^u
g(Z^\delta(\sigma(s)))r_j(K^\delta(\sigma(s))
V_j\left(\frac{L(\sigma,s)}{\delta}\right)ds
\end{eqnarray}
and
\begin{eqnarray}\label{d5122}
III_2=\frac{1}{{\delta}^{3/2}}\int_{t+\delta^{1-\gamma}}^u\int_0^1
g(Z^\delta(\sigma(s)))r_j(K^\delta(\sigma(s))
V_{jk}\left(\frac{R(v,\sigma,s)}{\delta}\right)
(X_k^\delta(s)-L_k(\sigma,s))  dvds.
\end{eqnarray}
The expectation of the first term on the event $[\tau_\delta>T]$
is small by the mixing lemma because the points $X^\delta(\sigma)$
and $L(\sigma,s)$ are at distance of order $\delta^{1-\gamma}$,
$\E[V_j(x)]=0$, and the mixing coefficient is rapidly decaying. On
the other hand, ${\mathbb P}[\tau_\delta<T]$ is small according to
Lemma \ref{lem-tau}. We conclude that
\begin{equation}\label{III-1-est}
\left|\E\left\{III_1\tilde\zeta\right\}\right|\le
C\delta\|G\|_1\left(\|\zeta\|_\infty+\E[\tilde\zeta]\right)(u-t).
\end{equation}
In order to estimate $III_2$ we write
\begin{eqnarray*}
V_{jk}\left(\frac{R(v,\sigma,s)}{\delta}\right)=
V_{jk}\left(\frac{L(\sigma,s)}{\delta}\right)+
\frac{1}{\delta}\int_0^v
V_{jkm}\left(\frac{R(\theta,\sigma,s)}{\delta}\right)
(X_m^\delta(s)-L_m(\sigma,s))d\theta,
\end{eqnarray*}
which splits $III_2=III_{21}+III_{22}$ (note that the $v$-variable
integrates out in $III_{21}$):
\begin{eqnarray}\label{d511}
III_{21}=\frac{1}{{\delta}^{3/2}}\int_{t+\delta^{1-\gamma}}^u
g(Z^\delta(\sigma(s)))r_j(K^\delta(\sigma(s))
V_{jk}\left(\frac{L(\sigma,s)}{\delta}\right)
(X_k^\delta(s)-L_k(\sigma,s))   ds,
\end{eqnarray}
and
\begin{eqnarray}\label{d512}
&&III_{22}=\frac{1}{{\delta}^{5/2}}\int_{t+\delta^{1-\gamma}}^u\int_0^1\int_0^v
g(Z^\delta(\sigma(s)))r_j(K^\delta(\sigma(s))
V_{jkm}\left(\frac{R(\theta,\sigma,s)}{\delta}\right)\\
&&~~~~~~~~~~~~~~\times(X_m^\delta(s)-L_m(\sigma,s))
(X_k^\delta(s)-L_k(\sigma,s))d\theta dvds.\nonumber
\end{eqnarray}
We  have linearized sufficiently to make $III_{22}$ be
deterministically small because of (\ref{L-X-close}):
\begin{equation}\label{III22}
|III_{22}|\le C\delta^{-5/2}\delta^{3-2\gamma}\|G\|_{2}(u-t).
\end{equation}
This leaves us with $III_{21}$ to take care of. This we do with
the help of (\ref{L-X-approx})
\begin{eqnarray*}
&&\!\!\!\!\!\!\!\!\!\!\!\!\!\!\!\!\!\!\!\!\!\!\!\!\!\!
X_i^\delta(s)-L_i(\sigma,s)=-\int_\sigma^s[K(v)-K(\sigma)]dv=-
\frac{1}{\sqrt{\delta}}\int_\sigma^s\int_\sigma^{v}
V_j\left(\frac{X^\delta(\rho)}{\delta}\right)d\rho dv\\
&& = -\frac{1}{\sqrt{\delta}}\int_\sigma^s(s-\rho)
V_j\left(\frac{X^\delta(\rho)}{\delta}\right)d\rho=-
\frac{1}{\sqrt{\delta}}\int_\sigma^s(s-\rho)
V_j\left(\frac{L(\sigma,\rho)}{\delta}\right)d\rho\\
&&-\frac{1}{\sqrt{\delta}}\int_\sigma^s(s-\rho)
\left[V_j\left(\frac{X^\delta(\rho)}{\delta}\right)-
V_j\left(\frac{L(\sigma,\rho)}{\delta}\right)\right]d\rho=
-\frac{1}{\sqrt{\delta}}\int_\sigma^s(s-\rho)
V_j\left(\frac{L(\sigma,\rho)}{\delta}\right)d\rho\\
&&-\frac{1}{\delta^{3/2}}\int_\sigma^s\int_0^1(s-\rho)
V_{jm}\left(\frac{R(v,\sigma,\rho)}{\delta}\right)
\left(X_m^\delta(\rho)-L_m(\sigma,\rho)\right)dvd\rho.
\end{eqnarray*}
that allows us to decompose $III_{21}=III_{21}^1+III_{21}^2$ with
\begin{eqnarray}\label{d7211}
III_{21}^1=-\frac{1}{{\delta}^{2}}\int_{t+\delta^{1-\gamma}}^u\int_\sigma^s
(s-\rho)g(Z^\delta(\sigma(s)))r_j(K^\delta(\sigma(s))
V_{jk}\left(\frac{L(\sigma,s)}{\delta}\right)
V_j\left(\frac{L(\sigma,\rho)}{\delta}\right)d\rho ds,
\end{eqnarray}
and
\begin{eqnarray}
&&III_{21}^2=-\frac{1}{{\delta}^{3}}\int_{t+\delta^{1-\gamma}}^u
\int_\sigma^s\int_0^1(s-\rho)
g(Z^\delta(\sigma(s)))r_j(K^\delta(\sigma(s))
V_{jk}\left(\frac{L(\sigma,s)}{\delta}\right)
V_{jm}\left(\frac{R(v,\sigma,\rho)}{\delta}\right)\nonumber\\
&&~~~~~~~~~~~~\times\left(X_m^\delta(\rho)-L_m(\sigma,\rho)\right)dvd\rho
ds.\label{d7112}
\end{eqnarray}
Again, $III_{21}^2$ is deterministically small because of
(\ref{L-X-close}):
\[
|III_{21}^2|\le
C\delta^{-3}\delta^{2-2\gamma}\delta^{3/2-\gamma}\|G\|_{2}(u-t)\le
C\delta^{1/2-3\gamma}\|G\|_{2}(u-t).
\]
Now, for $III_{21}^1$ we first compute
\begin{eqnarray*}
&&\E\{V_{jk}(x+y)V_j(y)\}=\frac{\partial^2}{\partial x_j\partial
x_k}\E\{V (x+y)V_j(y)\}=-\frac{\partial^2}{\partial x_j\partial
x_k}\E\{V_j (x+y)V(y)\}\\
&&~~~~~~~~~~~~~~~~~~~~~~~~~~=-\frac{\partial^3}{\partial
x_j^2\partial x_k}\E\{V(x+y)V(y)\}=-\Delta R_k^{VV}(x)
\end{eqnarray*}
and use the mixing lemma to conclude that
\begin{eqnarray*}
&&\left|\E\left\{III_{21}^1\tilde\zeta \right\}-
\frac{1}{{\delta}^{2}}\int_{t+\delta^{1-\gamma}}^u\int_\sigma^s
(s-\rho)g(Z^\delta(\sigma(s)))r_j(K^\delta(\sigma(s)) \Delta
R_{j}^{VV}\left(\frac{L(\sigma,s)-L(\sigma,\rho)}{\delta}\right)
 d\rho ds\right|\\
 &&\le C\delta\|G\|_1\E[\tilde\zeta](u-t).
\end{eqnarray*}
As before we change variables in the inner integrals above and
replace the argument $\sigma(s)$ of smooth functions appearing
above by $s$, as well replacing the limits of integration by their
values in the limit $\delta\to 0$, to conclude that
\begin{eqnarray*}
&&\left|\E\left\{III_{21}^1\tilde\zeta \right\}-
 \int_{t }^u\int_0^\infty
\theta g(Z^\delta(s))r_j(K^\delta(s)) \Delta
R_{j}^{VV}\left(-\theta
K^\delta(s)\right) d\theta ds\right|\\
 &&\le C\delta\|G\|_{4}\E[\tilde\zeta](u-t).
\end{eqnarray*}
Overall, the work of this section implies that
\begin{eqnarray}\label{III-est}
&&\left|\E\left\{III \tilde\zeta \right\}+
 \int_{t }^u\int_0^\infty
\theta g(Z^\delta(s))r_j(K^\delta(s)) \Delta
R_{j}^{VV}\left(\theta
K^\delta(s)\right) d\theta ds\right|\\
 &&\le C\delta^\alpha\|G\|_4\left(\|\zeta\|_\infty+\E[\tilde\zeta]\right)(u-t)
 \nonumber
\end{eqnarray}
which is nothing but (\ref{III-est-lem}).

\subsection{Estimate of $IV$}

We are now down to the last term $IV$ in (\ref{d22}) that is
estimated as $II$ with the help of the mixing lemma and no
additional expansions:
\begin{eqnarray}
&&IV=\frac{1}{{\delta}}\int\limits_{t+\delta^{1-\gamma}}^u\int\limits_{\sigma(s)}^s
g'(Z^\delta(\rho))r_j(K^\delta(\rho))S\left(\frac{X^\delta(\rho)}{\delta}\right)
V_j\left(\frac{X^\delta(s)}{\delta}\right)d\rho ds\label{e6}\\
&&~~~~+
\frac{1}{{\delta}}\int\limits_{t+\delta^{1-\gamma}}^u\int\limits_{\sigma(s)}^s
g(Z^\delta(\rho))r_{jm}(K^\delta(\rho))V_m\left(\frac{X^\delta(\rho)}{\delta}\right)
V_j\left(\frac{X^\delta(s)}{\delta}\right)d\rho ds.\nonumber
\end{eqnarray}
First, we observe using the mixing lemma and smoothness of $G$
that (compare to (\ref{II-II'}))
\begin{equation}\label{IV-IV'}
\left|\E\left\{IV\tilde\zeta\right\}-\E\left\{IV'\tilde\zeta\right\}\right|\le
C\delta^{\gamma_2}\|G\|_3\left[\|\tilde\zeta\|_\infty+\E\{\tilde\zeta\}\right](u-t)
\end{equation}
with
\begin{eqnarray}
&&IV'=\frac{1}{{\delta}}\int\limits_{t+\delta^{1-\gamma}}^u\int\limits_{\sigma(s)}^s
g'(Z^\delta(\sigma))r_j(K^\delta(\sigma))S\left(\frac{L(\sigma,\rho)}{\delta}\right)
V_j\left(\frac{L(\sigma,s)}{\delta}\right)d\rho ds\label{e6prime}\\
&&~~~~+
\frac{1}{{\delta}}\int\limits_{t+\delta^{1-\gamma}}^u\int\limits_{\sigma(s)}^s
g(Z^\delta(\sigma))r_{jm}(K^\delta(\sigma))V_m\left(\frac{L(\sigma,\rho)}{\delta}\right)
V_j\left(\frac{L(\sigma,s)}{\delta}\right)d\rho
ds=IV_1+IV_2.\nonumber
\end{eqnarray}
Now, we have by the mixing lemma
\begin{eqnarray*}
&&\left|\E\{IV_1\tilde\zeta\}-\frac{1}{{\delta}}
\int\limits_{t+\delta^{1-\gamma}}^u\int\limits_{\sigma(s)}^s
g'(Z^\delta(\sigma))r_j(K^\delta(\sigma))
R_j^{SV}\left(\frac{L(\sigma,s)-L(\sigma,\rho)}{\delta}\right)
 d\rho ds\tilde\zeta\right| \\
 &&\le
C\delta^\alpha\|G\|_2\left(\|\tilde\zeta\|_\infty+E\{\tilde\zeta\}\right)(u-t)
\end{eqnarray*}
and thus
\begin{equation}
\left|\E\{IV_1\tilde\zeta\}-
\int\limits_{t}^u\int\limits_{0}^\infty
g'(Z^\delta(s))r_j(K^\delta(s)) R_j^{SV}\left(-\theta K^\delta(s)
\right)
 d\theta ds\tilde\zeta\right|
\le
C\delta^\alpha\|G\|_2\left(\|\tilde\zeta\|_\infty+E\{\tilde\zeta\}\right)(u-t).\label{IV1-est}
\end{equation}
Similarly, for the other contribution we have
\begin{eqnarray*}
&&\left|\E\{IV_2\tilde\zeta\}+\frac{1}{{\delta}}\int\limits_{t+\delta^{1-\gamma}}^u\int\limits_{\sigma(s)}^s
g(Z^\delta(\sigma))r_{jm}(K^\delta(\sigma))R_{mj}^{VV}\left(\frac{L(\sigma,s)-L(\sigma,\rho)}{\delta}\right)
d\rho ds\tilde\zeta\right|\nonumber\\
 &&\le
C\delta^\alpha\|G\|_2\left(\|\tilde\zeta\|_\infty+E\{\tilde\zeta\}\right)(u-t)
\end{eqnarray*}
and
\begin{equation}
\left|\E\{IV_2\tilde\zeta\}+
\int\limits_{t}^u\int\limits_{0}^\infty
g(Z^\delta(s))r_{jm}(K^\delta(s))R_{mj}^{VV}\left(\theta
K^\delta(s) \right)
d\theta ds\tilde\zeta\right|
\le
C\delta^\alpha\|G\|_2\left(\|\tilde\zeta\|_\infty+E\{\tilde\zeta\}\right)(u-t).\label{IV2-est}
\end{equation}
Together, (\ref{IV1-est}) and (\ref{IV2-est}) imply
(\ref{IV-est-lem}). This finishes the proof of Lemma
\ref{lem-I-IV} and thus of Proposition~\ref{prop-mart} as well.
$\Box$


\end{document}